\documentclass[12pt]{amsart}
\usepackage{amsmath}
\usepackage{amssymb} % \mathfrak, \mathcal
\usepackage{amsthm}
\usepackage{mathtools}
\usepackage{mathrsfs} %\mathscr
\usepackage{dsfont} % \mathds
\usepackage{bm}
\usepackage{bbm}  % \mathbbm, \mathbbss, \mathbbmtt
\usepackage{cite}
\usepackage{enumerate}
\usepackage{units}
\usepackage{syntonly}
\usepackage{stdclsdv}
\usepackage{alltt}
\usepackage{fontenc}
\usepackage{fancyhdr}
\usepackage{eqnarray}
\usepackage{calc}
\numberwithin{equation}{subsection}
\theoremstyle{plain}
\newtheorem{theo}{Theorem}
\newtheorem*{theo*}{Theorem}

\theoremstyle{definition}

\newtheorem{lem}{Lemma}[section]

\newtheorem{prop}[lem]{Proposition}
\newtheorem{cor}[lem]{Corollary}

\theoremstyle{remark}

\newtheorem*{rem*}{Remark}
\newtheorem{rem}[lem]{Reamrk}

\begin{document}

\title{Locally Lipschitz graph property for lines}
\author{Xiaojun Cui}
%\date{December 10, 2013}

\address{Xiaojun Cui \endgraf
Department of Mathematics, Nanjing University,
Nanjing 210093, Jiangsu Province,  People's
Republic of China.}
\email{xcui@nju.edu.cn}

%\curraddr{}

\thanks{Supported by the National Natural Science Foundation of China
 (Grant 11271181), the Project Funded by the Priority Academic Program Development of Jiangsu Higher Education Institutions  (PAPD) and the Fundamental Research Funds for the Central Universities.}

 \abstract
On a non-compact, smooth, connected, boundaryless, complete Riemannian manifold $(M,g)$, one can define its ideal boundary by rays (or equivalently, Busemann functions). From the viewpoint of Mather theory, boundary elements could be regarded as  the static classes of Aubry sets,  and thus lines should be think as the semi-statics curves connecting different static classes. In Mather theory, one core property is Lipschitz graph property for Aubry sets and for some kind of semi-static curves. In this article, we prove a such kind of result for a set of lines which connect the same pair of boundary elements.
\endabstract
\maketitle

%\tableofcontents
\section*{Introduction}
Let $M$ be a smooth, non-compact, complete, boundaryless, connected  Riemannian manifold with Riemanian metric $g$. Let $TM$ be the tangent bundle and $\pi$ be the canonical projection of $TM$ onto $M$. The distance $d$ on $M$ is induced by the Riemannian metric $g$ and on $TM$ it is induced by the Sasaki metric $g_S$. We also use $L$ to denote the length of a curve with respect to $g$. Throughout this paper, all geodesic segments are always parametermized to be unit-speed.  By a ray, we mean a geodesic segment $\gamma:[0, +\infty) \rightarrow M$ such that $d(\gamma(t_1),\gamma(t_2))=|t_2-t_1|$ for any $t_1,t_2 \geq 0$. Throughout this paper, $| \cdot |$ mean Euclidean norms. To guarantee the existence of rays, $M$ must be non-compact. By definition, the Busemann function associated to a ray $\gamma$, is defined as
$$b_{\gamma}(x):=\lim_{t \rightarrow +\infty}[d(x, \gamma(t))-t].$$
Clearly, $b_{\gamma}$ is a Lipschitz function with Lipschitz constant 1, i.e.
$$|b_{\gamma}(x)-b_{\gamma}(y)| \leq d(x,y).$$
Moreover, it is proved in \cite{CC} that $b_{\gamma}$ is locally semi-concave with linear modulus (for the definition, see \cite{CS}).

By rays, or their Busemann functions, one could define the ideal boundary $M(\infty)$. By  definition, $M(\infty)$ is the set of  equivalent classes of rays, where two rays $\gamma_1$ and $\gamma_2$ are equivalent if and only if $b_{\gamma_1}=b_{\gamma_2}+ const.$. For more concrete information on $M(\infty)$, we refer to \cite{Ba}, \cite{Mar}. It should be noted that one can use horofunctions  or $dl$-functions, instead of rays (or equivalently, Busemann functions) to define other kinds of ideal boundary, for details, see \cite{Mar}, \cite{Ba} or \cite{CC}.

Now we introduce a more restrictive notation than ray. A geodesic $\gamma: \mathbb{R} \rightarrow M$ is called to be a line if
$$d(\gamma(t_1), \gamma(t_2))=|t_2-t_1|$$
for any $t_1, t_2 \in \mathbb{R}.$ On any non-compact complete Riemannian manifold $M$, for any point $x$ on it, there always exists at least one ray initiated from $x$. Comparing with rays, lines are much more rare  and there exist examples of non-compact complete Riemannian manifolds containing no line at all. Given a line $\gamma:\mathbb{R} \rightarrow M$, we can define the barrier function $B_{\gamma}$, under the motivation of Mather theory, to be
$$B_{\gamma}(x)=b_{\gamma^{+}}(x)+b_{\gamma^{-}}(x),$$
where $\gamma^{+},\gamma^{-}$ are two rays defined by
$\gamma^{\pm}:[0, +\infty) \rightarrow M, \gamma^{\pm}(t)=\gamma(\pm t)$. It should be mentioned that such kind of barrier functions has appeared in literature (e.g. [\cite{Pe}, Page 287, Line -15--Line-10],   [\cite{Sa}, Page 218, Line -3]) in the context of non-negative Ricci curvature. Here, we follow the ideas from \cite{CC}.

For any line $\gamma:\mathbb{R} \rightarrow M$, $- \gamma$ is defined by $(-\gamma)(t):=\gamma(-t)$,  for any $t \in \mathbb{R}$.

In Mather theory, the core property is the Lipschitz graph property of Aubry Sets (For simplicity, we only consider the autonomous version). In fact, it is easy to generalize  the Lipschitz graph  property by Mather's curve shortening lemma [\cite{Mat1}, Page 186, Lemma]  to the set:
$$\dot{N}_{A_1,A_2}:=\begin{Bmatrix}\cup (\gamma(t), \dot{\gamma}(t)): (\gamma(t), \dot{\gamma}(t))  \text{ are semi-static orbits} \\
\text{ with } \alpha\text{-limit set in } A_1 \text{ and } \omega\text{-limit set in }A_2
\end{Bmatrix} ,\hfill{(*)} $$
here $A_1$ and $A_2$ are any two fixed static classes. For terminologies or  more details on Mather's theory, we refer to \cite{Mat1}, \cite{Mat2}, \cite{Fa}.

In \cite{CC}, the authors made a simple, but  interesting observation that on a noncompact Riemannian manifold, one could regard the elements in $M(\infty)$ as the analogous of static classes in Mather's theory. This analogue is reasonable: For example, by weak KAM theory, every static class will determine a (globally defined) viscosity solution up to a constant; for elements in $M(\infty)$, it is also true. Thus, one would ask whether there are still some Lipschitz graph property for some invariant (with respect to geodesic flow) sets? This is the main motivation of this paper and we could answer this problem in a confirmed way.

Since the elements in $M(\infty)$ are represented by rays, thus for any line $\gamma$, we could think $\gamma$ as a geodesic connecting two rays: $\gamma|[\tau_1, \infty)$ and $(-\gamma)|[\tau_2, \infty)$ for $\tau_1,\tau_2 \in \mathbb{R}$.  Since $\gamma$ is a line, $\gamma|[\tau_1, \infty)$ and $\gamma|[\tau_2, \infty)$ determine a same element in $M(\infty)$ for any $\tau_1, \tau_2 \in \mathbb{R}$, the same is true for $-\gamma$. So the choice of $\tau_1, \tau_2$ is not crucial at all, we may choose $\gamma^-:= (-\gamma)|[0,\infty)$ and $\gamma^+:= \gamma|[0,\infty)$ as the  representations of the two elements in $M(\infty)$ which are connected by the line $\gamma$.

Analogous to $\dot{N}_{A_1,A_2}$ defined in $(*)$, we could define
$$\dot{N}_{\gamma^-,\gamma^+}:=\begin{Bmatrix}\cup (\xi(t), \dot{\xi}(t)): \xi(t)  \text{ are lines } \\
\text{ and }  \xi(t) \text{ connects } \gamma^- \text{ and }\gamma^+
\end{Bmatrix} \hfill{(**)}$$
 and let ${N}_{\gamma^-,\gamma^+}$ be the projection of $\dot{N}_{\gamma^-,\gamma^+}$ into $M$.
In other words, the trajectories  $(\xi, \dot{\xi})$ in $\dot{N}_{\gamma^-,\gamma^+}$ must satisfy $b_{\xi^\pm}=b_{\gamma^\pm}+ const.$.

Our main result in this paper is
\begin{theo}
$\pi:\dot{N}_{\gamma^-,\gamma^+} \rightarrow N_{\gamma^-,\gamma^+}$ is a locally bi-Lipschitz map.
\end{theo}

The main tool for the proof of this theorem is still the curve shortening lemma \cite{Mat1}.
\begin{proof}
In [\cite{CC}, Proposition 5.1], the authors proved that a trajectory $(\xi, \dot{\xi})$ lies in $\dot{N}_{\gamma^-,\gamma^+}$ if and only if $\xi \sim \gamma$. Here, $\sim $ is an equivalence relation defined by $\xi \sim \gamma$ if and only if  $\xi \prec \gamma $ and $\gamma \prec \xi$, and the relation $\prec$ is defined further as follows.

We say $\xi \prec \gamma$ if :

$\bullet$ $B_{\gamma}(\xi(t)) \equiv 0$;

$\bullet$ For any $\tau \in \mathbb{R}$, $\xi|[\tau, \infty)$ is a coray to $\gamma^+$ and $(-\xi)|[\tau, \infty)$ is a coray to $\gamma^-$. Here, for definition of the coray, we refer to \cite{Bu}, and here we only mention a useful property (e.g. [\cite{BG}, Proposition 2.7]): A ray  $\gamma_1$  is  a coray to the ray $\gamma$ if and only if $b_{\gamma}(\gamma_1(t_1))-b_{\gamma}(\gamma_1(t_2))=t_2-t_1$ for any $t_1, t_2 \geq 0$.

For any compact subset $S$, let  $TS:=\{ \cup TM_x: x \in S\}$ . We will prove that there exists a constant $C$ (depend on the compact subset $K$) such that $$d(\pi^{-1}(x_1) \cap \dot{N}_{\gamma^-,\gamma^+}, \pi^{-1}(x_2)\cap \dot{N}_{\gamma^-,\gamma^+}) \leq C d(x_1,x_2)$$ for any $x_1,x_2 \in K \cap {N}_{\gamma^-,\gamma^+}.$ First, we choose two connected open sets $U_1,U_2$ with compact closures and such that $D_1(K) \subset  U_1 $ and  $D_1( \bar{U_1}) \subset U_2 $, here, for any compact set $S$, $D_1(S)$ means the closed ball $D_1(S):=\{x \in M:d(x,S) \leq 1\}$, and $\bar{U_1}$ means the closure of $U_1$. Then the Riemannian metric $g$, restricted on the open set $U_2$, is uniformly  positive definite.

Now we reformulate Mather's curve shortening lemma [\cite{Mat1}, Page 186, Lemma] (see also [\cite{Mat2}, Page 1361--Page 1362] for further explantation)  as follows:

 \begin{lem}There exist $\epsilon, \delta, C>0$ such that if $\alpha, \beta:[t_0- \epsilon, t_0+\epsilon] \rightarrow M$ are two geodesic segments in $U_1$ with $d(\alpha(t_0), \beta(t_0)) \leq \delta$ and $d((\alpha(t_0), \dot{\alpha}(t_0)), (\beta(t_0), \dot{\beta}(t_0))) \geq C d(\alpha(t_0), \beta(t_0))$, then there exist $C^1$ curves $c_1,c_2:[t_0-\epsilon,t_0+\epsilon] \rightarrow U_2$ such that $c_1(t_0-\epsilon)=\alpha(t_0-\epsilon), c_1(t_0+\epsilon)=\beta(t_0+\epsilon), c_2(t_0-\epsilon)=\beta(t_0-\epsilon), c_2(t_0+\epsilon)=\alpha(t_0+\epsilon)$, and
 $$4 \epsilon >L(c_1)+L(c_2).$$
\end{lem}
For the proof of this lemma, we refer to \cite{Mat1}. Despite of some slight modifications (just restrict Mather's result to the setting of geodesic flow), the original proof  still goes through. We remark that $\epsilon, \delta$ in the lemma could be taken arbitrarily small.

Based on this lemma, we could go on to prove our main result as follows. Assume that our theorem is not true.  Then for any $\delta >0, C>0$, there exist two lines $\xi$ and $\eta$ with $(\xi, \dot{\xi}), (\eta,\dot{\eta}) \in \dot{N}_{\gamma^-,\gamma^+}$,  a real number $t_0$ (if necessary, we operate a time translation on $\xi$ or on $\eta$), such that $\xi(t_0) \in K, \eta(t_0) \in K,$
$d(\xi(t_0), \eta(t_0)) \leq \delta$ and $d((\xi(t_0), \dot{\xi}(t_0)), (\eta(t_0), \dot{\eta}(t_0))) \geq C d(\xi(t_0), \eta(t_0))$. Choose $\epsilon$ sufficiently small, then $\xi|[t_0-\epsilon, t_0+\epsilon] \subset U_1, \eta|[t_0-\epsilon, t_0+\epsilon] \subset U_1$. Choose  $\delta, C$ suitably, such that the conditions of curve shortening lemma (Lemma 0.1) is satisfied.  By the curve shortening lemma, there exist two $C^1$ curves  $c_1,c_2:[t_0-\epsilon,t_0+\epsilon] \rightarrow U_2$ such that $c_1(t_0-\epsilon)=\xi(t_0-\epsilon), c_1(t_0+\epsilon)=\eta(t_0+\epsilon), c_2(t_0-\epsilon)=\eta(t_0-\epsilon), c_2(t_0+\epsilon)=\xi(t_0+\epsilon)$, and
 $$4\epsilon >L(c_1)+L(c_2).$$

Thus,
\begin{eqnarray*}&&B_{\gamma}(c_1(t_0)) \\
&=&b_{\gamma^+}(c_1(t_0))+b_{\gamma^-}
(c_1(t_0))\\
&\leq&  b_{\gamma^+}(\eta(t_0+\epsilon))+d(c_1(t_0), \eta(t_0+\epsilon))+b_{\gamma^-}(\xi(t_0-\epsilon))+d(c_1(t_0), \xi(t_0-\epsilon))\\
&\leq &b_{\gamma^+}(\eta(t_0+\epsilon))+ b_{\gamma^-}(\xi(t_0-\epsilon))+L(c_1)\\
&= & b_{\gamma^+}(\eta(t_0))-\epsilon+ b_{\gamma^-}(\xi(t_0))-\epsilon+L(c_1)\\
&=& b_{\gamma^+}(\eta(t_0))+ b_{\gamma^-}(\xi(t_0))-2\epsilon+L(c_1), \hfill{(\#)}
\end{eqnarray*}
where the first equality follows from the fact that Busemann functions are Lipschitz with Lipschitz constant 1; the second equality holds because for any $\tau \in \mathbb{R}$, $\xi|[\tau, \infty),\eta|[\tau, \infty)$ are corays to $\gamma^+$ and $(-\xi)|[\tau, \infty),(-\eta)|[\tau, \infty)$ are corays to $\gamma^-$.
\end{proof}

Analogously, we could get
$$B_{\gamma}(c_2(t_0)) \leq b_{\gamma^+}(\xi(t_0))+ b_{\gamma^-}(\eta(t_0))-2\epsilon+L(c_2).\hfill{(\#\#)}$$

Combing inequalities $(\#)$ and $(\#\#)$, we obtain
 \begin{eqnarray*}&&B_{\gamma}(c_1(t_0))+B_{\gamma}(c_2(t_0))\\
 &\leq& b_{\gamma^+}(\eta(t_0))+ b_{\gamma^-}(\xi(t_0))-2\epsilon+L(c_1)+b_{\gamma^+}(\xi(t_0))+ b_{\gamma^-}(\eta(t_0))-2\epsilon+L(c_2)\\
 &=& B_{\gamma}(\eta(t_0))+B_{\gamma}(\xi(t_0))-4\epsilon+L(c_1)+L(c_2)\\
 &<&0,
 \end{eqnarray*}
 where the last inequality holds because $B_{\gamma}\xi \equiv 0, B_{\gamma}\eta \equiv 0$ and $L(c_1)+L(c_2)< 4\epsilon$. But it contradicts  the fact that $B_{\gamma}$ is a non-negative function. The contradiction proves Theorem 1. \qed

By the procedure of the proof of Theorem 1,  in fact we could get a stronger result as follows. For any line $\gamma$, let 
$$\dot{G}_{\gamma}:=\{\cup(\xi(t), \dot{\xi}(t)):\xi \text{ are lines with } \xi \prec \gamma\}$$
 and $G_{\gamma}$ be the projection of $\dot{G}_{\gamma}$ into $M$. The stronger result is
 
 \begin{prop}
 $\pi: \dot{G}_{\gamma} \rightarrow G_{\gamma}$ is a locally bi-Lipschitz map.
 \end{prop}
 
\begin{rem}
The injectivity of $\pi|\dot{N}_{\gamma^-,\gamma^+}$ or $\pi|\dot{G}_{\gamma}$ also follows from the differentiability of $B_{\gamma}$ on ${G}_{\gamma}$ [\cite{CC}, Theorem 3. 2) ].
\end{rem}
\begin{cor}
If ${G}_{\gamma}=M$, the by Theorem 1, $M$ will be foliated by the curves $\xi$ where $\xi \prec \gamma$. Moreover, the foliation is locally Lipschitz. In this case, $b_{\gamma^+}$ (also, $b_{\gamma^-}$) is a locally $C^{1,1}$ viscosity solution (in fact, classical solution) to the Hamilton-Jacobi equation
$$| \nabla u|_g=1,$$
where $\nabla$ is the gradient, and $|\cdot |_g$ is the norm on $TM$, both determined by the Riemannian metric $g$.
\end{cor}

\end{document}